\documentclass{amsart}

\usepackage{amsthm}
\usepackage{amssymb}
\usepackage{amsmath}

\newcommand{\labell}[1] {\label{#1}}
\newtheorem {Theorem}   {Theorem} 

\numberwithin{Theorem}{section}

\theoremstyle{definition}

\theoremstyle{remark}
\newtheorem{Remark}[Theorem]{Remark}

\newtheorem {Corollary}[Theorem]{Corollary}  

\expandafter\chardef\csname pre amssym.def at\endcsname=\the\catcode`\@ 
\catcode`\@=11 
\def\undefine#1{\let#1\undefined} 
\def\newsymbol#1#2#3#4#5{\let\next@\relax 
 \ifnum#2=\@ne\let\next@\msafam@\else 
 \ifnum#2=\tw@\let\next@\msbfam@\fi\fi 
 \mathchardef#1="#3\next@#4#5}
\def\mathhexbox@#1#2#3{\relax 
 \ifmmode\mathpalette{}{\m@th\mathchar"#1#2#3}% 
 \else\leavevmode\hbox{$\m@th\mathchar"#1#2#3$}\fi} 
\def\hexnumber@#1{\ifcase#1 0\or 1\or 2\or 3\or 4\or 5\or 6\or 7\or 8\or 
 9\or A\or B\or C\or D\or E\or F\fi} 
 
\font\teneufm=eufm10 
\font\seveneufm=eufm7
\font\fiveeufm=eufm5 
\newfam\eufmfam
\textfont\eufmfam=\teneufm 
\scriptfont\eufmfam=\seveneufm 
\scriptscriptfont\eufmfam=\fiveeufm 
 
\catcode`\@=\csname pre amssym.def at\endcsname 

\def	\eps	{\epsilon}

\def	\C {{\mathbb C}}
\def	\reals	{{\mathbb R}}
\def	\R	{{\mathbb R}}

\def	\T	{{\mathbb T}}
\def	\CP	{{\mathbb C}{\mathbb P}}

\def	\codim	{\operatorname{codim}}
\def	\rk	{\operatorname{rk}}

\def	\pr	{\operatorname{pr}}
\def	\SB	{\operatorname{SB}}
\def	\CL	{\operatorname{CL}}

\begin{document}

%%%%%%%%%%%%%%%%%%%%%%%%%%%%%%
%   TEXT FORMATTING

%%%%%%%%%%%%%%%%%%%%%%%%%%

%%%%%%%%%%%%%%%%%%%%%%%%%%

%%%%%%%%%%%           BEGINNING OF  TEXT

%%%%%%%%%%%%%%%%%%%%%%%%%%

\title[Periodic Orbits in Magnetic Fields]{Periodic Orbits in Magnetic Fields
in Dimensions Greater Than Two}

\author[Viktor L. Ginzburg]{Viktor L. Ginzburg}
\author[Ely Kerman]{Ely Kerman}
\address{Department of Mathematics, UC Santa Cruz, 
Santa Cruz, CA 95064}
\email{ginzburg@math.ucsc.edu; ely@cats.ucsc.edu}

\date{January, 1999}

\thanks{The work is partially supported by the NSF and by the faculty
research funds of the University of California, Santa Cruz.}

\subjclass{Primary: 58F05, 58F22}

\bigskip

\begin{abstract}
The Hamiltonian flow of the standard metric Hamiltonian with respect to 
the twisted symplectic structure on the cotangent bundle describes the
motion of a charged particle on the base. We prove that under certain 
natural hypotheses the number of periodic orbits on low energy levels for 
this flow is at least the sum of Betti numbers of the base.

The problem is closely related to the existence question for periodic
orbits on energy levels of a proper Hamiltonian near a Morse--Bott 
non-degenerate minimum. In this case, when some extra requirements are met,
we also give a lower bound for the number of periodic orbits. 
Both of these questions are very similar to the Weinstein conjecture 
but differ from it in that the energy levels may fail to have contact 
type.

We give a simple proof of the fact that bounded sets in the cotangent bundle to
the torus, with a twisted symplectic structure, have finite Hofer--Zehnder
capacity. As a consequence, we obtain the existence of periodic orbits 
on almost all energy levels for magnetic fields on tori. 
\end{abstract}

\maketitle

\section{Introduction} 
The motion of a charge in a magnetic field on a manifold is given 
by the Hamiltonian vector
field of the standard metric Hamiltonian on the cotangent bundle 
with respect to a non-standard symplectic structure. This so-called
twisted symplectic structure is the sum of the standard one and the
pull-back of the magnetic field two--form. 

One of the objectives of the present paper is to prove the existence of 
periodic orbits on
low energy levels for such Hamiltonian systems, provided that the metric
and the magnetic field two--form satisfy a certain condition.
The condition is that the magnetic field two--form is symplectic and
compatible with the metric. For example, this requirement is met
when the manifold, with these structures, is K\"{a}hler. More specifically, 
we show (Theorem 
\ref{cor:magnetic}) that under this condition the number of periodic 
orbits on every low energy level, when the orbits are non-degenerate, 
is no less than the sum of Betti numbers of the manifold.
In the degenerate case, a similar lower bound is obtained in 
\cite{kerman} in terms of the minimal number of critical points of
a function on the manifold.

The problem of existence of periodic orbits for 
symplectic magnetic fields can be generalized as follows. 
Consider a proper function on a symplectic manifold and assume that 
this function has a Morse-Bott non-degenerate minimum along a symplectic 
submanifold. Then the problem is to find a lower bound for
the number of periodic orbits of the Hamiltonian flow 
on the levels near the minimum. For example, when the submanifold
is just a point the answer is given by Weinstein's theorem,
\cite{weinstein-1973}. We provide a lower bound (Theorem 
\ref{thm:min}) when the orbits are non-degenerate and again a certain 
compatibility condition is satisfied. The case of degenerate orbits
is treated in \cite{kerman} using a method relying on Moser's
proof, \cite{Moser-1976,Bottkol}, of Weinstein's theorem. 

Both of these problems are similar to the Weinstein conjecture, 
\cite{we:conj}, on the existence of periodic orbits on contact type
hypersurfaces. (See, e.g., \cite{fhv,hv,ho-ze:book,MS,vi:Theorem,
vi:functors} for more information and further references.) 
However, the essential difference is that the energy levels
in question may fail to have contact type. For example, in the case
of the magnetic field, the twisted symplectic form is never exact
on the energy level when the magnetic field is not exact and
the dimension of the base is greater than two. This fact makes
the problem difficult to solve by standard symplectic topology techniques.

The problem of existence of periodic orbits on almost all levels
is accessible by making use of symplectic capacities.
For example, if one can show that bounded sets in the ambient manifold
have finite Hofer--Zehnder capacity, it follows that a proper Hamiltonian
has periodic orbits on almost all energy levels, 
\cite[Section 4.2]{ho-ze:book}. In Section \ref{sec:tori}, we prove that
the capacity of bounded sets is finite for a twisted symplectic form on 
the cotangent bundle to a torus for any closed magnetic field two--form
(Theorem \ref{thm:capacity}). This implies the ``almost existence'' 
of periodic orbits for magnetic fields on tori (Corollary 
\ref{cor:torus}).

This paper is one of very few (see also \cite{bialy,BT,kerman,Lu,polt})
focusing on magnetic fields on manifolds of dimension greater
than two. Much more is known about the existence of periodic orbits for
magnetic fields on surfaces. The reader interested in the review of
these results from the symplectic topology perspective should
consult \cite{gi:Cambr}. The necessary symplectic geometry 
material can be found in \cite{ho-ze:book,MS}.

\section{Periodic orbits on low energy levels} 

\subsection{Periodic orbits near a minimum.}
\labell{sec:min}
Let $(W, \omega)$ be a symplectic manifold and let
$H\colon W\to \R$ be a smooth function. Assume that $H$ has 
a Bott non-degenerate minimum along a compact symplectic
submanifold $M\subset W$. The normal bundle $\nu$ to $M$
in $W$ is a symplectic vector bundle of dimension $2m=\codim M$.
The Hessian
$d^2H\colon \nu\to\R$ is a positive--definite fiberwise
quadratic form on $\nu$. Thus at every point $x\in M$ we have
$m$ canonically defined eigenvalues of the Hessian $d^2_xH$ on the fiber
of $\nu$ over $x$ with respect to the linear symplectic form 
$\omega^\perp_x$ on this fiber.  The eigenvalues of 
$d^2_x H$ are equal at every point $x\in M$ if and only if
the linear Hamiltonian flow of $d^2_xH$ with respect to 
$\omega^\perp_x$ is periodic with all its orbits having the same period. 
In this case, the flow gives rise to a free $S^1$-action on 
$S\nu=\{d^2H=1\}$ and hence to the principal $S^1$-bundle
$\pr\colon S\nu\to S\nu/S^1$. 

Recall that a closed integral curve of a vector field (or a
line field)  is said to be non-degenerate
if its Poincar\'{e} return map does not have unit as an eigenvalue.
Denote by $\SB(N)$ the sum of Betti numbers of a manifold $N$.

\begin{Theorem}
\labell{thm:min}
Assume that for every $x\in M$ all eigenvalues of $d^2_xH$ with
respect to $\omega^\perp_x$ are equal. Then
for a sufficiently small $\eps>0$, the number of periodic orbits 
of the Hamiltonian flow of $H$ on the level $\{ H=\eps \}$
is at least $\SB( S\nu/S^1)$ if the $S^1$-bundle
$\pr$ is trivial and at least $\SB(S\nu)/2+1$ otherwise,
provided that all of the periodic orbits are  non-degenerate.
\end{Theorem}

The theorem will be proved in Section \ref{sec:proof}. Note that the bundle
$\pr$ is always non-trivial when $m>1$. The reason is that its restriction
to a fiber $\CP^{m-1}$ of the bundle $S\nu/S^1\to M$ is the Hopf
fibration.

\begin{Remark}
\labell{rmk:min}
The fiberwise Hamiltonian flow  of $d^2H$ can alternatively be described as the Hamiltonian 
flow of $d^2H$ with respect to the Poisson structure on the total space of 
$\nu$ given by the family of fiberwise symplectic forms $\omega^\perp_x$,
$x\in M$.

As follows from the results of \cite{kerman}, when the periodic 
orbits are not required
to be non-degenerate, the number of periodic orbits is 
still greater than or equal to $\CL(S\nu/S^1)=\CL(M)+m$. Here 
$\CL$ stands for the 
cup--length and the action of $S^1$ is given, after suitable
rescaling, by the linear Hamiltonian flow of $d^2H$, described above.
For $W=\R^{2n}$ and $M$ a point, this is a particular case of Weinstein's
theorem \cite{weinstein-1973}; see also \cite{Moser-1976}.

It seems to be likely that the bound of Theorem \ref{thm:min}
can be significantly improved. (For example, as stated, this bound is not
necessarily higher than the bound from below for the degenerate case 
mentioned above, \cite{kerman}.) We conjecture that under the hypothesis
of the theorem, the number of periodic orbits is at least
$\SB(S\nu/S^1)=m\SB(M)$. (The equality follows from \cite[Theorem 2.5,
p. 233]{Hu}.)

\end{Remark}

\subsection{Periodic trajectories of a charge in a magnetic field.}
\labell{sec:magn} 
Let $M$ be a Riemannian manifold and let $\sigma$ be a symplectic 
form on $M$ (a magnetic field). Then
$\omega=d\lambda+\pi^*\sigma$ is symplectic on  $W=T^*M$. Here
$d\lambda$ is the standard symplectic form on $T^*M$ and 
the map $\pi\colon T^*M\to M$ is the natural projection. Take
$H\colon T^*M\to \R$ to be the standard metric Hamiltonian.
More explicitly, let us identify $TM$ and $T^*M$ by means of the
Riemannian metric on $M$. Then $H(X)=g(X,X)/2$, where $g$ is the
metric.
The Hamiltonian flow of $H$ with respect to $\omega$ describes
the motion of a charge on $M$ in the magnetic field $\sigma$.
(The reader interested in more details should consult, e.g.,
\cite{gi:Cambr}.) We say that a symplectic form $\sigma$ and a metric
$g$ are compatible if there exists an almost complex structure
$J$ such that $g(X,Y)=\sigma(X,JY)$ for all $X$ and $Y$ and $J$ is
$g$-orthogonal (cf., \cite[Section 4.1]{MS}). This condition is 
equivalent to that at every point all eigenvalues of $g$ with 
respect to $\sigma$ are equal. 
In Section \ref{sec:proof} we will prove the following

\begin{Theorem}
\labell{cor:magnetic}
Assume that $\sigma$ is compatible with the Riemannian metric
on $M$, e.g., $M$ is K\"{a}hler with these structures. For a sufficiently 
small $\eps>0$, the 
number of periodic orbits on the energy level $\{H=\eps\}$ is at least
$\SB(M)$, provided that all of the orbits are non-degenerate.
\end{Theorem}

\begin{Remark}
\labell{rmk:magnetic}
Similarly to Theorem \ref{thm:min}, the lower bound given
by Theorem \ref{cor:magnetic} 
is probably far from sharp. In the non-degenerate case
there should be conjecturally at least $\SB(STM/S^1)=m\SB(M)$ periodic orbits
on every low energy level, where $2m=\dim M$. When $M$ is a surface, this
is exactly the statement of the Theorem \ref{cor:magnetic}, which in this
case was originally proved in \cite{gi:FA}. Note also that 
under the hypothesis of Theorem \ref{cor:magnetic},
there are at least $\SB(M)+1$ periodic orbits, when
$\chi(M)=0$ but $M\neq \T^2$. This can be seen easily from the proof of
the theorem.

If the periodic orbits are not assumed to be non-degenerate,
the lower bound of Theorem \ref{cor:magnetic} should be
replaced by $\CL(M)+m$, \cite{kerman}.

When the energy value $\eps>0$ is not small, the level 
$\{H=\eps\}$ may fail to carry a periodic orbit.
See \cite{gi:MathZ} for $m=1$ and \cite[Example 4.2]{gi:review}
for $m\geq 2$.
\end{Remark}

\subsection{The proofs of theorems \ref{thm:min} and \ref{cor:magnetic}}
\labell{sec:proof}
The theorems will follow from a more general result, proved in 
\cite{gi:FA}, which we now state.

Let $E$ be a compact odd--dimensional manifold with a free circle action.
Denote the quotient $E/S^1$ by $B$ and the natural projection
$E\to B$ by $\pr$. Recall that a two--form $\eta$ on $E$ is said
to be maximally non-degenerate if it has a one--dimensional null--space,
denoted in what follows by $\ker \eta_x$, at every point $x\in E$.
Equivalently, this means that $\eta$ has maximal possible rank (equal
to $\dim E-1$) at every point of $E$. Recall also that a characteristic
of $\eta$ is an integral curve of the line field $\ker\eta$.

\begin{Theorem}
\labell{thm:general}
Let $\eta$ be a closed maximally non-degenerate two--form on $E$
such that the cohomology class $[\eta]$ lies in the image
$\pr^*(H^2(B))$. 
Assume that the field of directions $\ker\eta$ is $C^1$-close to
the fibers of $\pr$ and that all of the closed characteristics of $\eta$
are non-degenerate. Then $\eta$ has at least $\SB(B)$ closed 
characteristics if the principle $S^1$-bundle $\pr$ is trivial
and at least $\SB(E)/2+1$ closed characteristics otherwise.
\end{Theorem}

\begin{Remark}
\labell{rmk:general}
We conjecture that in the non-degenerate case, regardless of whether
$\pr$ is trivial or not,
the number of closed characteristics is greater than or equal to 
$\SB(B)$, which would give a much higher lower bound
than that of Theorem \ref{thm:general}. 
Note that this conjecture would imply the hypothetical lower bound
$m\SB(M)$ mentioned in Remarks \ref{rmk:min} and \ref{rmk:magnetic}.
The conjecture is proved in \cite{gi:FA} for a surface $B$ and
in \cite{gi:MathZ} for the case where
$\eta$ is a $C^0$-small perturbation of the pull--back of
a symplectic form on $B$. The latter result
indicates that the requirement that $\ker\eta$ is $C^1$-close
to the fibers can perhaps be replaced by $C^0$-closeness. 

Note in this connection that the above conjecture
is erroneously claimed to be proved by the first of the authors in 
\cite{gi:MathZ} (Corollary 3.8). In fact, Corollary 3.8 does not follow 
from Theorem  2.7 of that paper for the reasons outlined in Remark 
\ref{rmk:perturb} below.  At present, the conjecture (and hence 
Corollary 3.8 of \cite{gi:MathZ}) appears to be an open problem.
\end{Remark}

\begin{Remark}
\labell{rmk:perturb}
Under the hypothesis of Theorem \ref{thm:general}, the ``horizontal
component'' of the average of $\eta$ is the pull-back of a symplectic 
form on $B$. If $\eta$ were close to this pull-back, 
Theorem \ref{thm:general}, with an improved lower bound, would follow from
\cite[Theorem 2.7]{gi:MathZ}. However, the form $\eta$ may fail to
be close to this pull-back and as a consequence 
\cite[Theorem 2.7]{gi:MathZ} does not apply. In fact, the main point of Theorem \ref{thm:general} is that the form 
$\eta$ is not assumed to be close to the pull--back of any two-form on 
$B$. In other words, one may think of $\eta$ as a ``Hamiltonian 
perturbation'' of a non-Hamiltonian vector field generating the 
$S^1$-action. This is exactly what makes Theorem \ref{thm:general}, in 
contrast with the results of \cite{gi:MathZ},  applicable to magnetic 
flows in higher dimensions. 
\end{Remark}

\begin{Remark}
It is easy to see that the cohomology condition on $\eta$ cannot
be omitted -- without this condition $\eta$ may fail to have closed
characteristics. (See \cite{gi:FA} for a more detailed discussion.)

If the closed characteristics of $\eta$ are not required to be 
non-degenerate, the low bound of Theorem \ref{thm:general}
should be replaced by the minimal number of critical points of
a smooth function on $B$, \cite{kerman}. In both of these lower 
bounds only the integral curves close to fibers and winding only once along 
the fibers are counted. In contrast with the lower bound of Theorem 
\ref{thm:general},
the lower bound of \cite{kerman} is apparently sharp for the number of
periodic orbits in the class specified.
\end{Remark}

\begin{proof}[Proof of Theorem \ref{thm:general}]
When $B$ is a surface, the theorem is proved in \cite{gi:FA}. Hence, in
what follows we will assume that $2k+1=\dim E\geq 5$. 

Furthermore, as is shown in \cite{gi:FA}, the number of closed 
characteristics is at least $\SB(E)/2$. Since
$\SB(E)/2=\SB(B)$ when the bundle $\pr$ is trivial, we will assume 
from now on that $\pr$ is not trivial.

We need to recall some details of the argument from \cite{gi:FA}.
Consider closed characteristics of $\eta$ which are close to the orbits of the
$S^1$-action and wind only once along the orbits. 
There exists a Morse--Bott function $f\colon E\to \R$ such that
the closed characteristics of $\eta$ of this type comprise exactly the 
critical manifolds of $f$. Hence it is sufficient to prove
that $f$ has at least $\SB(E)/2+1$ critical manifolds.

Denote by $b_i$ the Betti numbers of $E$ and by $\mu_i$ the number 
of critical manifolds of $f$ of index $i$. 
All (co)homology groups will be taken with real coefficients and
we will use the convention that $\mu_i=0$ whenever $i$ is negative. 

Since the critical manifolds of $f$ are 
circles, the Morse--Bott inequalities for $f$ turn into
\begin{equation}
\labell{eq:Morse-Bott}
\mu_i+\mu_{i-1}\geq b_i \text{ for } i=0,\ldots, 2k+1 .
\end{equation}

The inequalities \eqref{eq:Morse-Bott}
can be refined when $i=1$ and $i=2k$. Namely, we 
claim that
\begin{equation}
\labell{eq:Morse-Bott2}
\mu_1\geq b_1 \text{ and } \mu_{2k-1}\geq b_{2k},
\end{equation}
Let us prove the first of the inequalities \eqref{eq:Morse-Bott2}
(cf., \cite{gi:FA}).
Denote by $L\subset H_1(E)$ the subspace generated by the critical 
manifolds of index zero and $L'\subset H_1(E)$ a complementary subspace 
to $L$, i.e., $L\oplus L'=H_1(E)$. The dimension of this space does
not exceed the number of critical circles of index one:
\begin{equation}
\labell{eq:proj}
\mu_1\geq \dim L'
. 
\end{equation}
The critical manifolds of 
$f$ are close to the orbits of the $S^1$-action. As a consequence, the 
projections of critical manifolds to $B$ are contained in small balls. 
In particular, these projections are contractible and $\pr_*(L)=0$ in 
$H_1(B)$. Thus \eqref{eq:proj} implies that 
$\mu_1\geq\dim L'\geq\dim\pr_*(H_1(E))$.
However, $\pr_*(H_1(E))=H_1(B)$ and so 
$$
\mu_1\geq \dim H_1(B).
$$
Finally, since $\pr$ is a non-trivial bundle, $b_1=\dim H_1(B)$. Hence
$\mu_1\geq b_1$ which proves the first inequality
in \eqref{eq:Morse-Bott2}. The second inequality follows from the first
one with $f$ replaced by $-f$.

Adding up the Morse--Bott inequalities \eqref{eq:Morse-Bott}
for all $i=0,\ldots,2k+1$ except $i=1$ and $i=2k$ and the refined 
Morse--Bott inequalities \eqref{eq:Morse-Bott2}, we obtain 
$$
2\sum_{i=0}^{2k}\mu_i -(\mu_0+\mu_{2k})\geq \sum_{i=0}^{2k+1}b_i
.
$$
As a consequence,
$$
\sum\mu_i\geq \SB(M)+\frac{1}{2}(\mu_0+\mu_{2k}).
$$
Since $f$ must have a local minimum  and a local maximum,
the second term in the latter formula is
greater than or equal to one. Therefore, $\sum\mu_i\geq \SB(E)/2+1$ which 
completes the proof.
\end{proof}

\begin{proof}[Proof of Theorem \ref{thm:min}.] 
Throughout the proof we will keep the notations and conventions of
Section \ref{sec:min}. 
Using the symplectic neighborhood theorem, let us identify a neighborhood of 
$M$ in $W$ with a neighborhood $U$ of the zero section in the total 
space of the normal bundle
$\nu$ so that the symplectic structure on $U$ is linear and equal to
$\omega^\perp$ on the fibers of $\nu$.
Without loss of generality we may assume that $U$ contains the level 
$E=\{ d^2H=1\}$. 

Denote by $\varphi_\eps$ the fiberwise dilation 
$y\mapsto \eps y$ in the fibers on $\nu$. Let $E_\eps$ be the level
$\varphi_\eps^*H/\eps^2=1$ and $\psi_\eps\colon E\to E_\eps$ be the
fiberwise central projection. Consider the vector field $X_\eps$ on 
$E$ obtained by pushing forward the Hamiltonian vector field of 
$H/\eps^2$ on the level $\{H=\eps^2\}$ to $E_\eps$ by means 
of $(\psi_\eps\varphi_\eps)^{-1}$. 

Clearly, the assertion of the theorem is equivalent to 
$X_\eps$ having the required number of periodic orbits.

It is not hard to show (see Remark \ref{rmk:limit} below)
that $X_\eps\stackrel{C^1}{\to} X_0$ as $\eps \to 0$, where
$X_0$ is the fiberwise Hamiltonian vector field
of $d^2H$ with respect $\omega^\perp$. In other words, $X_0$
is the Hamiltonian vector field of $d^2H$ for the Poisson structure which
has the fibers of $\nu$ as its leaves and is given by $\omega^\perp$.
Furthermore, the vector field $X_\eps$ generates the null-space 
line field of the two--form $\eta_\eps=(\psi_\eps\varphi_\eps)^*\omega$.

By the hypothesis on the eigenvalues of $d^2H$, the flow of $X_0$
gives rise to a free $S^1$-action on $E$.
The form $\eta_\eps$ satisfies the cohomology condition of Theorem 
\ref{thm:general}. Thus, by Theorem \ref{thm:general}, when
$\eps>0$ is small enough the form $\eta_\eps$
has at least $\SB(E/S^1)$ or $\SB(E)/2+1$ closed characteristics,
depending on whether $\pr\colon E\to E/S^1$ is trivial or not and provided 
that the closed characteristics are non-degenerate. This completes the proof
of Theorem \ref{thm:min}.
\end{proof}

\begin{Remark}
\labell{rmk:limit}
Set $H_\eps=\varphi_\eps^*H/\eps^2$. Then 
$H_\eps\stackrel{C^k}{\to} d^2 H$ for any $k\geq 0$  as $\eps\to 0$ 
on a small neighborhood $U$ of the zero section. Furthermore,
consider the family of symplectic forms 
$\omega_\eps=\varphi^*_\eps\omega$. This family of forms does not
converge as $\eps\to 0$, but the corresponding Poisson structures
do. To be more specific, denote by $\omega_\eps^{-1}$ the Poisson
structure on $U$ corresponding to the symplectic structure $\omega_\eps$
and by $\omega_0^{-1}$ the Poisson structure on $U$ whose leaves
are the fibers of $\nu$ and which is equal to $\omega^\perp$ on the 
fibers. One can verify (see \cite{kerman}) that
$\omega_\eps^{-1}\stackrel{C^k}{\to}\omega_0^{-1}$ on $U$ as
$\eps\to 0$ for any $k\geq 0$. Therefore, the Hamiltonian vector
field of $H_\eps$ with respect to $\omega_\eps$ $C^k$-converges to
the Hamiltonian vector field of $d^2H$ with respect to $\omega_0^{-1}$.
\end{Remark}

\begin{proof}[Proof of Theorem \ref{cor:magnetic}.] 
The normal bundle $\nu$ to $M$ is canonically isomorphic to $T^*M$.
Furthermore, let us identify the tangent and cotangent bundles to $M$ 
by means of $\sigma$. Then $\omega^\perp_x=\sigma_x$ for all $x\in M$.
The compatibility condition of the corollary is equivalent to
that for every $x\in M$ the eigenvalues of the metric on $T_xM$
with respect to $\sigma_x$ are equal. Thus the hypothesis 
of Theorem \ref{thm:min} holds under the compatibility condition.

As a consequence, Theorem \ref{thm:min} gives a lower bound for the number
of periodic orbits on the energy level $\{H=\eps\}$ in terms of
the sum or Betti numbers of $E=ST^*M$ or $E/S^1$. 

The $S^1$-bundle $\pr\colon E\to E/S^1$ is trivial only when $M=\T^2$. 
In this case $\SB(E/S^1)=\SB(M)$. Thus let us assume that $\pr$ is non-trivial.

If $\chi(M)=0$, we have $\SB(E)/2=\SB(M)$ and Theorem \ref{cor:magnetic}
follows from Theorem \ref{thm:min} with an even higher lower bound.

Assume that $\chi(M)\neq 0$. We claim that 
\begin{equation}
\labell{eq:sb1}
\SB(E)/2+1\geq \SB(M)
.
\end{equation}
Denote by $b_i$ the Betti numbers of $E$ and by $\beta_i$ the Betti
numbers of $M$. The condition $\chi(M)\neq 0$ implies that
\begin{equation*}
%\labell{eq:Betti}
b_i=
\begin{cases}
\beta_i\text{ if } 0\leq i\leq n-1\\
\beta_{i-n+1}\text{ if } n\leq i\leq 2n-1,
\end{cases}
\end{equation*}
where $n=\dim M$. Adding up these equalities for all $i=0,\ldots, 2n-1$,
we obtain
\begin{equation}
\labell{eq:sb2}
\SB(E)\geq 2\SB(M)-(\beta_0+\beta_n)
.
\end{equation}
Without loss of generality we can assume that $M$ is connected. Then
$\beta_0=\beta_n=1$ and \eqref{eq:sb1} follows from \eqref{eq:sb2}.
This completes the proof of Theorem \ref{cor:magnetic}.
\end{proof}

\section{Magnetic fields on tori}
\labell{sec:tori}

In this section we give a simple proof of the fact that 
for a charge on a torus there are periodic orbits
on almost all energy levels.
Let $\sigma$ be a closed two-form on a torus $\T^n$. 
(Note that $n$ can be odd.) As in Section
\ref{sec:magn}, consider the twisted symplectic structure 
$\omega=d\lambda+\pi^*\sigma$ on $W=T^*\T^n$.

\begin{Theorem}
\labell{thm:capacity}
Every bounded set in $W$ has finite Hofer--Zehnder capacity.
\end{Theorem}

\begin{Remark}
The reader interested in the definition and properties
of the Hofer--Zehnder capacity should consult \cite{ho-ze:book}. 
Theorem \ref{thm:capacity} is not new. In fact, the theorem
has been known to experts for quite some time. When 
the cohomology class $[\sigma]$ is rational, the theorem follows from
\cite[Theorem 1.2]{jiang2}. For $\sigma$ 
symplectic, the theorem is proved in \cite[Lemma 5.3]{gi:Cambr}.
When $\sigma$ is the pull-back of a two-form under a projection
$\T^n\to\T^{n-1}$, the theorem becomes a particular case of a result of 
G. Lu, \cite[Theorem E]{Lu}. 
Finally, for $[\sigma]\neq 0$, Theorem \ref{thm:capacity} follows from
\cite[Theorem C]{Lu}. 
\end{Remark}

\begin{proof}Note first that to prove the theorem for $\sigma$ it 
suffices to prove the theorem for any form in the cohomology class 
$[\sigma]$. Indeed, the fiberwise shift 
by a one-form $\alpha$ sends bounded sets to bounded sets and transforms 
the twisted symplectic form $\omega$ into the form
$\omega-\pi^*d\alpha$. 

Hence without loss of generality, we may assume that $\sigma$ is 
a translation--invariant form $\T^n$. In other 
words, in some coordinates $x_1,\ldots, x_n$ on $\T^n$, we have
\begin{equation}
\labell{eq:linear}
\sigma=\sum a_{ij}dx_i\wedge dx_j ,
\end{equation}
for some constants $a_{ij}$. 

For an exact form $\sigma$, i.e., when $a_{ij}=0$, the theorem is
proved in \cite{jiang1}. This fact is very easy to see,
\cite[Proposition 4, p. 136]{ho-ze:book}: A bounded set
in $T^*S^1$ is contained in an annulus and the latter can be embedded
into a disc. By taking the product, we conclude that a bounded set in
$T^*\T^n$ can be symplectically embedded into a polydisc. A polydisc
has finite capacity and, as a consequence of monotonicity, a bounded 
set in $T^*\T^n$ has finite capacity.

Thus we may assume that $\sigma$ given by \eqref{eq:linear} is non-zero.

We claim that $(W, \omega)$ is symplectomorphic to
the product $\R^{2k}\times W_1$ with $k\geq 1$,
where $\R^{2k}$ is equipped with the standard symplectic form and
$W_1=\R^{n-2k}\times  \T^n$ is given a translation--invariant symplectic 
form.

Let us prove the claim. Consider the universal covering
$\tilde{W}=\R^{2n}$ of $W$ with the pull-back linear symplectic form 
$\tilde{\omega}$. Let $L$ be the inverse image of $\T^n$ in $\tilde{W}$
and let $L^\perp$ be the symplectic orthogonal complement to $L$
with respect to $\tilde{\omega}$. Pick a linear subspace $E$
in $L^\perp$ which is transversal to $L^\perp \cap L$.
The space $E$ is symplectic because the null--space of 
$\tilde{\omega}|_{L^\perp}$ is exactly $L^\perp\cap L$. Moreover,
$\dim E> 0$. This follows from the fact that $\rk \tilde{\omega}|_L >0$, and so
$L$ is not Lagrangian and $L^\perp\neq L$. Fix a symplectic subspace
$\tilde{W}_1$ in $\tilde{W}$ which contains $L$ and is transversal 
to $E$. Thus $\tilde{W}$ is symplectomorphic to $E\times \tilde{W_1}$

The decomposition $\tilde{W}=E\times \tilde{W_1}$ induces
the required direct product decomposition of $W$. To see this, note that 
$W=\tilde{W}/\Gamma$, where $\Gamma$ is a discrete subgroup in $L$.
Thus $W=E\times W_1$, where $W_1=\tilde{W}_1/\Gamma$.
It is clear that $E$ is symplectomorphic to $(\R^{2k}, \omega_0)$ 
with $2k>0$ and that the resulting symplectic structure on $W_1$ is
translation--invariant.

Observe now that every bounded subset of $W_1$ can be symplectically
embedded into $\T^{2(n-k)}$ with some translation--invariant symplectic
structure. As a result, every bounded open set in $W$ is
symplectomorphic to a bounded open set in $\R^{2k}\times \T^{2(n-k)}$,
where $2k>0$. The latter open sets have finite capacity as proved in
\cite{fhv} and \cite{Ma}; see also \cite{Lu} for further generalizations. 
(It is essential that $2k>0$ and hence the space $E$
is non-trivial: the torus $\T^{2(n-k)}$ may have infinite capacity, 
\cite[Section 4.5]{ho-ze:book}.)
\end{proof}

As in Section \ref{sec:magn}, let $H\colon T^*\T^n\to\R$ be a
metric Hamiltonian. Finiteness of capacity implies (see
\cite[Section 4.2, Theorem 4]{ho-ze:book}) ``almost existence'' of 
periodic orbits:

\begin{Corollary}
\labell{cor:torus}
Almost all, in the sense of measure theory, levels of $H$ carry a
periodic orbit.
\end{Corollary}

To put this general benchmark result in perspective, let us state
some recently proved more subtle theorems concerning the existence of
periodic orbits in the magnetic problem on higher--dimensional 
tori. (See, e.g., \cite{gi:Cambr} for a review of results for
surfaces.)

According to the first result, due to Polterovich, \cite{polt},
\emph{when $\sigma\neq 0$, there exists a sequence of 
positive energy values $c_k\to 0$ such that on every level $\{ H=c_k \}$ 
there exists a {\rm contractible} periodic orbit.} This is a rather deep
theorem: Corollary \ref{cor:torus} guarantees the existence of periodic
orbits on almost all levels of $H$, but does not guarantee that these
orbits are contractible. For instance, if $\sigma=0$ periodic orbits
are just closed geodesics and a metric on $\T^n$ can easily fail to have
contractible closed geodesic (e.g., a flat metric). This result also holds
for any compact manifold $M$ with $\chi(M)=0$ in place of $\T^n$.
Note also that as follows from a theorem of G. Lu, \cite[Theorem C]{Lu},
almost all levels of $H$ carry a contractible periodic orbit when
$[\sigma]\neq 0$.

The second theorem, a version of Hopf's rigidity, is due to Bialy, 
\cite{bialy}. By Bialy's theorem, \emph{every energy
level of $H$ carries an orbit with conjugate points, provided that
the metric is conformally flat and again $\sigma\neq 0$.} This fact
is related to the question of existence of contractible periodic orbits
because every such orbit (with non-zero Maslov index) would have 
conjugate points. Thus Bialy's theorem serves as indirect evidence in
favor of the affirmative answer to the existence question.

In conclusion note that when $\sigma$ is exact, all high energy levels have
contact type and thus carry periodic orbits, \cite{hv}. (See
\cite[Remark 2.3]{gi:Cambr} for more details.)

\end{document}